\pdfoutput=1

\documentclass[reqno]{amsart}

\usepackage[utf8]{inputenc}
\usepackage[T1]{fontenc}
\usepackage[USenglish]{babel}

\usepackage{braket}
\usepackage[per-mode=symbol]{siunitx}
\usepackage{mathtools}
\usepackage{booktabs}
\usepackage[foot]{amsaddr}

\usepackage{csquotes}
\usepackage[inline]{enumitem}

\usepackage{graphicx}
\usepackage{xcolor}

\usepackage[color=yellow!,textsize=small]{todonotes}

\mathtoolsset{showonlyrefs}

\newtheorem{thm}{Theorem}

\newtheorem{algo}[thm]{Algorithm}

\newcommand{\Natural}{\mathbb{N}}
\newcommand{\Real}{\mathbb{R}}
\newcommand{\Sphere}[1]{\mathbb{S}^{#1}}
\newcommand{\Probability}{\mathbb{P}}
\newcommand{\Expectation}[2]{\mathbb{E}_{#1}\left[ #2 \right]}

\newcommand{\upL}{\mathrm{L}}
\newcommand{\Lebesgue}[2]{\upL^{#1}(#2)}

\newcommand{\Dictionary}{\mathcal{D}}
\newcommand{\Hilbert}{\mathcal{H}}
\newcommand{\Transpose}{\mathrm{T}}

\renewcommand{\to}{\rightarrow}
\newcommand{\Differential}[1]{\:\mathrm{d}#1}

\DeclareMathOperator{\Span}{span}
\DeclareMathOperator*{\Argmax}{argmax}

\DeclarePairedDelimiter{\pParantheses}{(}{)}

\newcommand{\B}[1]{\pParantheses*{#1}}
\DeclarePairedDelimiter{\pAbs}{\lvert}{\rvert}
\newcommand{\Abs}[1]{\pAbs*{#1}}
\DeclarePairedDelimiter{\pNorm}{\lVert}{\rVert}
\newcommand{\Norm}[1]{\pNorm*{#1}}

\DeclarePairedDelimiter{\pInnerProduct}{\langle}{\rangle}
\newcommand{\InnerProduct}[2]{\pInnerProduct*{#1,#2}}

\DeclarePairedDelimiter{\pIntCC}{[}{]}
\newcommand{\IntCC}[2]{\pIntCC*{#1,#2}}
\DeclarePairedDelimiter{\pIntCO}{[}{)}
    \newcommand{\IntCO}[2]{\pIntCO*{#1,#2}}

\DeclarePairedDelimiter{\pIntOO}{(}{)}
\newcommand{\IntOO}[2]{\pIntOO*{#1,#2}}

\renewcommand{\epsilon}{\varepsilon}
\newcommand{\eps}{\epsilon}
\renewcommand{\phi}{\varphi}
\renewcommand{\theta}{\vartheta}
\renewcommand{\rho}{\varrho}
\renewcommand{\subset}{\subseteq}

\newcommand{\Closure}[1]{\overline{#1}}
\newcommand{\AbelPoisson}{Q}
\newcommand{\AbelPoissonK}[2]{\AbelPoisson^{(#1)}_{#2}}
\newcommand{\SurfaceMeas}{\omega}

\makeatletter
\newcommand{\ie}{i.\,e\@ifnextchar{.}{}{.\@}}
\newcommand{\iid}{i.\,i.\,d\@ifnextchar{.}{}{.\@}}
\newcommand{\eg}{e.\,g\@ifnextchar{.}{}{.\@}}
\newcommand{\wrt}{w.\,r.\,t\@ifnextchar{.}{}{.\@}}
\newcommand{\cf}{cf\@ifnextchar{.}{}{.\@}}
\makeatother

\begin{document}
\title[Three-dimensional simulation of nonwoven fabrics]{Three-dimensional simulation of nonwoven fabrics using a greedy approximation of the distribution of fiber directions}

\author{Simone Gramsch}
\address[SG]{Department Transport Processes\\Fraunhofer ITWM\\Fraunhofer-Platz 1\\67663 Kaiserslautern\\Germany}
\email{simone.gramsch@itwm.fraunhofer.de}
\author{Max Kontak}
\email{kontak@mathematik.uni-siegen.de}
\author{Volker Michel}
\address[MK,VM]{Geomathematics Group\\Department of Mathematics\\University of Siegen\\Walter-Flex-Str. 3\\57068 Siegen\\Germany}
\email{michel@mathematik.uni-siegen.de}
\thanks{MK and VM gratefully acknowledge the financial support by the School of Science and Technology of the University of Siegen, Germany.}


\begin{abstract}
    An elementary algorithm is used to simulate the industrial production of a fiber of a 3-dimensional nonwoven fabric.
    The algorithm simulates the fiber as a polyline where the direction of each segment is stochastically drawn based on a given probability density function (PDF) on the unit sphere.
    This PDF is obtained from data of directions of fiber fragments which originate from computer tomography scans of a real non-woven fabric.
    However, the simulation algorithm requires numerous evaluations of the PDF.
    Since the established technique of a kernel density estimator leads to very high computational costs, a novel greedy algorithm for estimating a sparse representation of the PDF is introduced.
    Numerical tests for a synthetic and a real example are presented.
    In a realistic scenario, the introduced sparsity ansatz leads to a reduction of the computation time for 100 fibers from nearly 40 days to 41 minutes.
\end{abstract}

\maketitle

\section{Introduction}
According to EDANA, the European Disposables and Nonwovens Association, a nonwoven is defined as \textquote[{cited in \cite{Mao:2010}}]{\ldots a sheet of fibres, continuous filaments, or chopped yarns of any nature or origin, that have been formed into a web by any means, and bonded together by any means, with the exception of weaving or knitting.}.
Their impact on our daily live is large due to their versatile properties.
Nonwovens can be absorbent, antistatic, breathable, conductive or non-conductive, elastic, flame resistant, impermeable or permeable, smooth, and stiff to name but a few.
These properties are achieved by combining raw materials with specific production processes.

\begin{figure}
    \centering
    \includegraphics[width=0.75\textwidth]{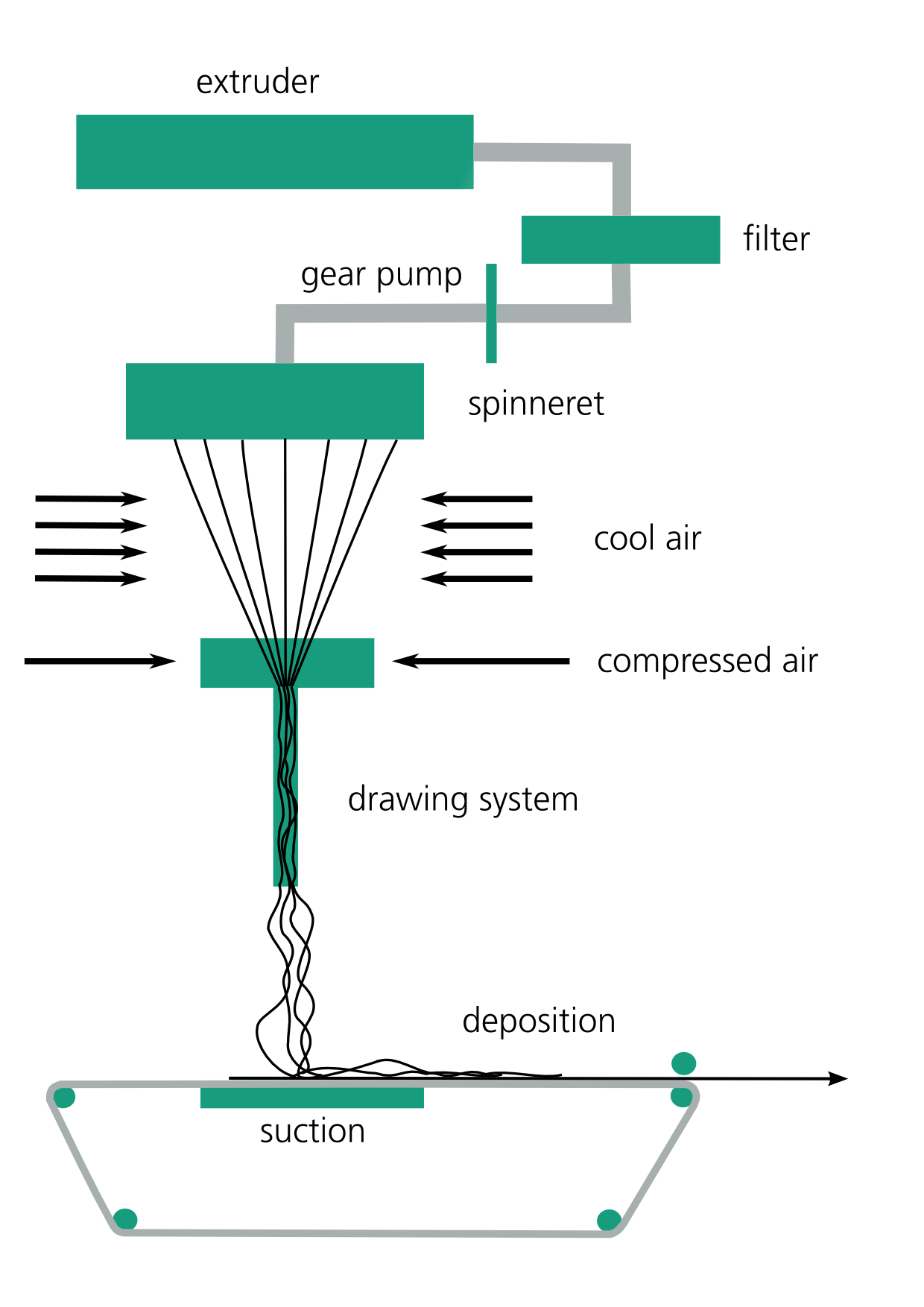}
    \caption{Sketch of a nonwoven production process (spunbond)} \label{f:process}
\end{figure}

There are three main nonwoven manufacturing processes:
dry-lay processes, wet-lay processes, and extrusion processes.
In the following, we concentrate on a typical extrusion process, the so-called spunbond process.
A sketch of this process can be found in Figure~\ref{f:process}.
In this production process, a polymer melt is extruded through spinnerets.
The fibers evolving from the spinnerets are cooled and stretched by air.
Additionally, they are swirled around by turbulent air streams until they are deposited on a perforated conveyor belt.
By suction, they are fixed on the belt to form a random web.
In order to produce a nonwoven fabric, additional process steps like bonding and finishing have to be implemented.
For further details we refer to \cite{Albrecht:2003}.

The advantage of nonwovens to be specially designed to have specific properties also has a drawback.
The market of nonwovens demands more and more customized products, so the development cycles are shortened accordingly.
Hence, simulating nonwoven production processes is a mathematical key technology that enables engineers to design the processes with respect to customer-specific needs.

Several papers have already dealt with the mathematical modeling of this specific production process.
One model that is based on the description of the physical process itself can be found in \cite{Klar:2009,Wegener:2015}, which is also used in the simulation software FIDYST (fiber dynamics simulation tool, cf.\ \cite{Gramsch:2015}).
In principle, one has to deal with a two-way coupling problem of the aerodynamic forces and the fiber dynamics, which is not solvable to industrial scales due to the required resolution of the mesh (cf.\ \cite{Marheineke:2011}).
Therefore, the model only incorporates a one-way coupling by using an air drag model.
Up to now, this approach neglects the interaction of one filament with itself or other filaments on the conveyor belt that have already been laid down.
Therefore, the simulated laydown of filaments on the conveyor belt is two-dimensional.

Due to the high computational effort needed for simulating this model, a real nonwoven consisting of thousands of fibers cannot be simulated in a reasonable amount of time.
Hence, in \cite{Goetz:2007,Klar:2009} a different approach was developed.
The two-dimensional fiber laydown is described by a surrogate model based on a stochastic differential equation.
The parameters that arise in this equation can be estimated from the FIDYST simulation of one single fiber.
Then, a considerable number of fibers representing the microstructure of a real nonwoven can be computed in a short time using the stochastic differential equation model. 

On the one hand, a two-dimensional simulated fabric still represents typical nonwovens, which are very flat, in an excellent way.
On the other hand, several properties of nonwovens originate in the fact that the fibers are lying on top of each other, which cannot be achieved by a two-dimensional simulation.
Thus, a three-dimensional model is desirable.

For that reason, in \cite{Klar:2012} the concept of the stochastic surrogate model was extended to a three-dimensional setting.
Although possessing good theoretical properties, this model lacks the possibility of estimating the parameters directly using a FIDYST simulation, since the output of the simulation is only two-dimensional.
Parameter estimation for the three-dimensional model by combining two-dimensional FIDYST data and three-dimensional data obtained by CT scans has been presented in \cite{Grothaus:2014} under the restriction that a non-moving belt or negligible small belt speeds are considered.

\begin{figure}
        \centering
        \input{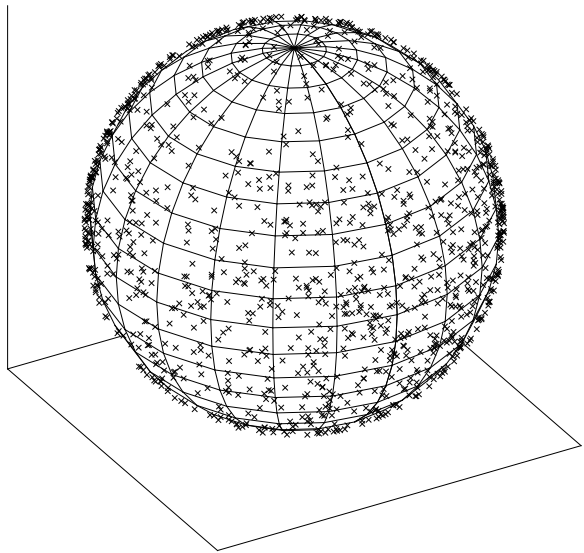}
        \caption{\num{2000} out of \num{9600558} points from a data set of directions in a real nonwoven fabric}
        \label{f:ct:pts}
\end{figure}

In section~\ref{s:math} of this paper, we will present an algorithm for the three-dimensional simulation of nonwovens based on the analysis of fiber directions inside a real nonwoven.
Our algorithm generates a polyline by sampling from a probability density function, which has to be estimated first.
This estimation is based on a data set generated by image analysis techniques.
Here, we use data obtained by 3D-microtomography with a resolution of \SIrange{1}{70}{\micro\metre} voxel side length and a maximal size of the sample of \SIrange{1}{10}{\cubic\milli\metre}.
A reduced version of the data set generated by the computer tomography (CT) scan is shown in Figure~\ref{f:ct:pts}.
Note that in a post-processing step, the raw CT data are translated to filament directions.
Since directions in the three-dimensional space can be equivalently described as points on the sphere, every data point on the sphere in Figure~\ref{f:ct:pts} corresponds to the filament direction in one single voxel.
For a more detailed description of the data set that is used, consult section~\ref{s:num} of this paper.
The resulting number of occurring voxels is extremely large due to the high resolution.
Hence, efficient strategies for analyzing the data must be developed.

In section~\ref{s:kde}, we will point out that, in principle, a standard technique of nonparametric statistics for the estimation of probability density functions, namely kernel density estimators, can be used to implement the simulation algorithm.
However, it is not efficient to use these estimators for very large data sets as in our case.
In particular, the fiber simulation requires numerous evaluations of the estimated density.

Therefore, we will derive a greedy algorithm for the estimation of probability density functions in section~\ref{s:greedy}.

In section~\ref{s:num}, we will first study the convergence of the greedy algorithm empirically and will then show that the efficiency of our simulation algorithm can be enhanced extremely by using the greedy algorithm to estimate the distribution of fiber directions.
We fill furthermore present first results about the validation of fiber laydown simulation methods like FIDYST, which can be performed using the presented greedy algorithm.

\section{Mathematical setting and the algorithm}\label{s:math}
No matter if mathematically simulated fibers are considered or if CT-based data of nonwoven fabrics are available, in each case we are confronted with a set of $N\in\Natural$ directions in the plane (2D case) or the space (3D case).
We will model these data $X_1,\ldots,X_N$ as independent and identically distributed realizations of an $\Sphere{d-1}$-valued random variable $X$, where $d=2$ or $d=3$.
By $\Sphere{d-1} \coloneqq \Set{ \xi\in\Real^d | \Abs{\xi} = 1}$, we denote the unit sphere in $\Real^d$.
Moreover, the distribution of $X$ shall be absolutely continuous with respect to the surface measure $\SurfaceMeas$ on $\Sphere{d-1}$, that is, there exists a probability density function (PDF) $f:\Sphere{d-1}\to\IntCO{0}{\infty}$, such that for every subset $A\subset\Sphere{d-1}$ we have
\begin{align}
    \Probability\B{X\in A} & = \int_A f\B{\xi} \Differential{\SurfaceMeas}\B{\xi},
\end{align}
where $\Probability\B{X\in A}$ denotes the probability that the random variable $X$ attains values in the subset $A$.
Note that we do not impose further assumptions on the PDF, especially not on any parametrization, such that the estimation of the PDF is a problem in the field of nonparametric statistics.

In the following, we propose a data-driven algorithm for the simulation of nonwoven fabrics.
The concept of the algorithm is to start at an arbitrary point in $\Real^d$ and to choose the next point of a fiber by sampling a direction from the estimated PDF and moving into this direction by a given step width.
By construction, the resulting fiber approximately has the same distribution of directions as it occurs inside the given data set used for estimating the PDF.
This idea translates to an algorithm as follows.
\begin{algo}\label{a:sim}
    Let an estimate $\hat{f}$ of the PDF $f$ as well as a discretization parameter $s\in\IntOO{0}{\infty}$ be given.

    Generate a discretization $\B{Z_j}_{j=0,1,\ldots} \subset \Real^d$ of a fiber by the following iteration:
    \begin{enumerate}
        \item Set $j\coloneqq 0$ and choose an initial point $Z_0\in\Real^d$.
        \item Sample a direction $Y_{j+1} \in \Sphere{d-1}$ from the estimated PDF $\hat{f}$.\label{a:sim:2}
        \item Set $Z_{j+1} \coloneqq Z_{j} + sY_{j+1}$.
        \item If the desired number of points is reached: stop.\\
              Otherwise: increase $j$ by 1 and return to step~\ref{a:sim:2}.
    \end{enumerate}
\end{algo}
Note the resemblance of the algorithm to the stochastic time discrete approximation of stochastic differential equations, in particular the Euler-Maruyama method, where Gaussian pseudo-random numbers are used in a similar way if the stochastic differential equation involves Brownian motion (cf. \cite{Kloeden:1992}, chapter~9).
However, it is not trivial if there is an interpretation of our algorithm as an Euler-Maruyama scheme for a certain stochastic differential equation and even if it was, it is completely unclear what such a stochastic equation or a stochastic process solving this equation would look like.

Note that the approach for the fiber simulation presented here differs essentially from other methods such as the one described in \cite{Klar:2009}.
For example, we do not solve a given stochastic differential equation which includes an explicit term for the belt movement as it is done in other cases.
However, in our numerical experiments in section~\ref{s:num}, we will use a PDF which is determined from data of a real nonwoven and, thus, incorporates the non-uniform distribution of the fiber directions with a dominance of the direction of the belt movement.
In this respect, the belt movement implicitly becomes an ingredient in our Algorithm~\ref{a:sim}.

A standard tool to solve the problem of estimation of PDFs is a so-called kernel density estimator, which is briefly discussed in the following section.
\section{Analysis of directions with kernel density estimators}\label{s:kde}
Since their introduction in \cite{Parzen:1962}, kernel density estimators (KDEs) have been a tool that has often been used for the estimation of PDFs based on given data sets $X_1,\ldots,X_N$.
Being introduced on the real line, extensions to the multivariate case (see \cite{Cacoullos:1966}) were straight-forward, which made this method superior to older methods like, for example, the finite difference approximation in \cite{Rosenblatt:1956}.
In \cite{Hall:1987}, KDEs were transferred to the case of spherical data $X_1,\ldots,X_N\in\Sphere{d-1}$ for $d\geq 2$, which is of particular importance for the presented application and is, thus, the foundation of the following considerations.

In analogy to the latter reference, in this paper a \emph{kernel} is the member of a family of integrable functions $K_\lambda: \IntCC{-1}{1}\to\IntCO{0}{\infty}$, which depend on a so-called \emph{concentration parameter} $\lambda\in\IntOO{0}{\infty}$, and have the following properties, that have also been considered in \cite{Freeden:2002} to obtain approximate identities on the sphere:
\begin{enumerate}[label=(\alph*)]
    \item normalization:
          \begin{align}
              \int_{-1}^1 K_\lambda\B{t} \Differential{t} & = 1, \label{e:ker:norm}
          \end{align}
    \item concentration at $1$: for every $c\in\IntOO{-1}{1}$, we have
          \begin{align}
              \lim_{\lambda\to\infty} \int_{c}^1 K_\lambda\B{t} \Differential{t} & = 1.\label{e:ker:conc}
          \end{align}
\end{enumerate}
The latter is different from the definition of kernels on $\Real$, where a concentration at $0$ is required, because for $\xi,\eta\in\Sphere{d-1}$ the propositions $\Abs{\xi-\eta} = 0$ and $\xi\cdot\eta = 1$ are equivalent and the kernel will always be used in the form $K_\lambda\B{\xi\cdot\eta}$.
Here, by $x\cdot y$ we denote the Euclidean inner product of two vectors $x,y\in\Real^d$.
Under certain further assumptions on the family of kernels it was proved in \cite{Hall:1987} that the associated KDE, which is defined as
\begin{align}
    \hat{f}\B{\xi} & \coloneqq \frac{1}{N} \sum_{n=1}^N K_\lambda\B{X_n\cdot\xi},\label{e:kde}
\end{align}
is an asymptotically unbiased and consistent estimator of the PDF $f$ if $\lambda$ is chosen properly in dependence of $N$.

We use the $d$-dimensional Abel-Poisson kernel $\AbelPoissonK{d}{h}:\IntCC{-1}{1}\to\IntOO{0}{\infty}$,
\begin{align}
    \AbelPoissonK{d}{h}\B{t} & \coloneqq \frac{1}{\SurfaceMeas\B{\Sphere{d-1}}} \frac{1-h^2}{\B{1+h^2-2ht}^{d/2}} \label{e:ap:def}
\end{align}
for $d=3$ in the following to obtain a KDE for a given CT data set, which consists of points on $\Sphere{2}$.
The term $\SurfaceMeas\B{\Sphere{d-1}}$ denotes the surface of the unit sphere $\Sphere{d-1}$ in $\Real^d$, which can be explicitly computed as
\begin{align}
    \SurfaceMeas\B{\Sphere{d-1}} & = \frac{2\pi^{d/2}}{\Gamma\B{\frac{d}{2}}},
\end{align}
which reduces to the well-known results $\SurfaceMeas\B{\Sphere{1}} = 2\pi$ and $\SurfaceMeas\B{\Sphere{2}} = 4\pi$ in the two- and the three-dimensional case, respectively.

Note that the Abel-Poisson kernel fulfills the conditions \eqref{e:ker:norm} and \eqref{e:ker:conc} if we set
\begin{align}
    K_\lambda\B{t} & \coloneqq \AbelPoissonK{d}{1-1/\lambda}\B{t}.
\end{align}
For further properties of the Abel-Poisson kernel in the three-dimensional case, see for example \cite[Section~5.1]{Michel:2013}.

\begin{figure}[t]
    \centering
    \input{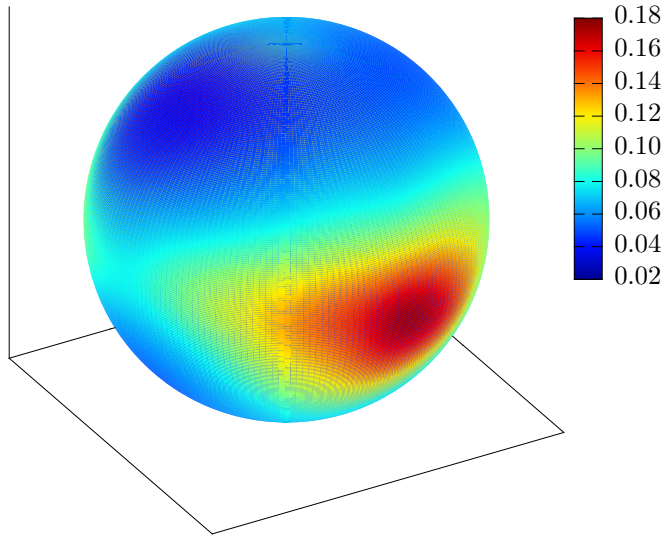}
    \caption{Plot of the KDE obtained from the data set where the Abel-Poisson kernel with parameter $h=0.9$ was used as the kernel in \protect\eqref{e:kde}.}\label{f:ct:kde}
\end{figure}

In Figure~\ref{f:ct:kde}, the KDE with respect to the CT data set is depicted on $\Sphere{2}$.
The parameter $h=0.9$ is more or less arbitrarily chosen, since a rigorous method for parameter selection is not known,
whereas heuristics, which are available, could of course be used to determine the parameter (for an overview on this subject, see \cite{Jones:1996} and the references therein).
The figure shows that there exists a dominant direction of the fibers, which corresponds to the direction in which the belt moves in the production process.

Concerning Algorithm~\ref{a:sim},
it turns out in numerical experiments that the use of a KDE as estimator of the PDF inside this method is not practicable from the computational point of view.
The crucial point is step~\ref{a:sim:2} of the algorithm, where sampling from the estimated PDF $\hat{f}$ is needed.

In our experiments---where we worked with a real CT data set on $\Sphere{2}$ with $N=\num{9600558}$ directions---we used the well-known acceptance-rejection method (see, \eg, \cite[Section~II.3]{Devroye:1986}) for sampling.
Simpler sampling methods like inversion sampling (see, \eg{}, \cite[Section~II.2]{Devroye:1986}) are unfortunately not applicable on spherical domains.

The acceptance-rejection method is based on the following idea (\cf{} \cite[Theorem~II.3.1]{Devroye:1986}):
If $\B{X,U}$ is a $\Sphere{d-1}\times\Real$-valued random variable, which is uniformly distributed on the set $A=\Set{\B{\eta,t} | \eta\in\Sphere{d-1}, 0\leq t \leq f\B{\eta}}$, then $X$ has the density $f$.
In practice, a sample $(\xi,u)$ of $(X,U)$ can be obtained by first generating a sample $(\tilde{\xi},\tilde{u})$ from a uniform distribution on $\tilde{A} = \Set{(\eta,t) | \eta\in\Sphere{d-1}, 0\leq t \leq C}$ for some upper bound $C$ of $f$.
If $\tilde{u}\leq f(\tilde{\xi})$, the sample is \emph{accepted} and $(\xi,u) := (\tilde{\xi},\tilde{u})$ is taken as a sample from the uniform distribution on $A$.
On the other hand, if $\tilde{u} > f(\tilde{\xi})$, the sample is \emph{rejected} and the generated pair is not taken into account for the samples from the uniform distribution on $A$.
Since methods for generating uniform distributions both on $\Sphere{d-1}$ and on the interval $[0,C]$ are well-known and implemented in well-established software like the GNU Scientific Library \cite{Galassi:2009}, samples from a uniform distribution on $\tilde{A} = \Sphere{d-1}\times [0,C]$ can easily be generated both in theory and in practice.

Note that the efficiency of the acceptance-rejection method depends heavily on the quality of the upper bound $C$ of the PDF.
When an estimation of the PDF is given by a KDE, a simple and rigorous method to obtain an upper bound is the application of the triangle inequality to the sum in \eqref{e:kde}.
Unfortunately, it turns out that this bound is by far too large and the computational effort is drastically high, since a vast majority of points is rejected in the acceptance-rejection method (\eg{}, \num{204} evaluations of the KDE were needed for one sample in numerical experiments).
In the numerical experiments presented in this paper, the upper bound was obtained by evaluation of the KDE on a equiangular grid with \num{10000} points on the sphere and multiplication of the result by \num{1.1} to bound maxima of the KDE, which may occur between the grid points.
This leads to an average of \num{2.4} evaluations of the KDE to obtain one sample in our numerical experiments, which improves the computational effort of the acceptance-rejection method in this case very much.

Note that for a single evaluation of the KDE in \eqref{e:kde}, the kernel has to be evaluated $N$ times, which affects the practicability of KDEs massively if $N$ is large.
Hence, since $N\approx \num{e7}$ is realistic in our application, the use of KDEs leads to very high computational costs.
In particular, in our numerical tests, approximately \num{2.4e7} evaluations of the Abel-Poisson kernel were needed to obtain one single sample $Y_j$ in step~\ref{a:sim:2} of Algorithm~\ref{a:sim}.
Additionally, by using the Abel-Poisson kernel, we already have a kernel at hand that admits a closed representation, whereas other kernels on the sphere (for examples, see \cite[Chapter~7]{Freeden:2009}) are given as a Legendre series and do not permit such a representation.
This increases the computational effort of this method even further, if different kernels are chosen.
Since a real nonwoven fabric can consist of hundreds of fibers and, in principle, the fibers can be arbitrarily long, we face a big problem in the practical realization of the method.
Specific numbers for computation times in comparison to the new method, which is developed in the following section, can be found in Table~\ref{t:times} later in this article.

To reduce the number of kernel evaluations, we consider the concept of sparsity as studied, for example, in the field of compressive sensing (\cf{} \cite[Chapter~2]{Foucart:2013}).
This means, in our case, that the aim is to estimate the PDF by a linear combination
\begin{align}
    \sum_{n=1}^N \alpha_n K_\lambda\B{X_n\cdot\xi}
\end{align}
of kernels, where most of the coefficients $\alpha_n$ are zero, or at least, very small.
Although the obtained result is technically not a kernel density estimator, we will call such a representation a \emph{sparse kernel density estimator}.
Note that this representation shows certain similarities to the approximation of functions in Sobolev spaces on the sphere with spherical splines (see \cite{Freeden:1981}).
In the literature, strategies to achieve a sparse representation of a function include optimization methods, greedy methods, and thresholding-based methods (\cf{} \cite[Chapter~3]{Foucart:2013}).
Sparse KDEs based on optimization methods were studied in \cite{Banerjee:2010,Chen:2013,Hong:2010}.
The arising optimization problems for the $N$ unknowns $\alpha_n$ may also be difficult to handle in practice if $N$ is large as in the presented application.
Greedy methods have only rarely been used in the estimation of PDFs so far:
a method for the parametric case with Gaussian kernels on $\Real$ is considered in \cite{Verbeek:2003}, which can obviously not be applied to the presented problem, since no parametrization of the PDF is assumed to be known.
The nonparametric case is the subject of \cite{Rosset:2002} with an approach based on the linearization and minimization of a so-called negative log-likelihood loss function.
Although the minimization of this specific loss function is a common method in statistics, apart from the linearization, it brings with it the difficulty that optimization methods have to be used in every iteration of the algorithm.
Furthermore, in that approach parts of the loss function have to be recalculated in every iteration, which demands a high computational effort for a large number of points.
These two difficulties do not arise with the new algorithm presented below, which can be interpreted as a minimization of an $\upL^2$-error.

\section{Analysis of directions by means of a greedy approximation}\label{s:greedy}
Greedy algorithms have been a successful tool in the sparse approximation of functions in Banach and Hilbert spaces (for a survey, see \cite{Temlyakov:2008}).
Several variants of these algorithms have also been used in signal processing as \emph{matching pursuits} (\cf{} \cite{Mallat:1993}) and in the solution of inverse problems as \emph{regularized functional matching pursuits} and its variants (\cf{} \cite{Fischer:2012,Michel:2014,Michel:2016}).
Note that greedy algorithms have also been used in statistics before, but usually for regression problems (\cf{} \cite{Barron:2008,Friedman:1981}).
In \cite{Friedman:1984}, density estimation in $\Real^d$ was considered, relying heavily on the Euclidean structure.

It is the aim of this section to transfer the basic concept of greedy algorithms to the estimation of PDFs.
The classical (or \emph{pure}) greedy algorithm for the approximation of functions in Hilbert spaces reads as follows (\cf{} \cite[Chapter~2]{Temlyakov:2008}).

\begin{algo}[Pure Greedy Algorithm, PGA]\label{a:pga}
    Let $\Hilbert$ be a Hilbert space with inner product $\InnerProduct{\cdot}{\cdot}$ and norm $\Norm{\cdot}$.
    Choose a so-called dictionary $\Dictionary\subset\Set{ g\in \Hilbert | \Norm{g}=1 }$.
    To approximate an element $f\in\Hilbert$, iteratively generate a sequence of approximations $(f_k)_{k\in\Natural_0}$ in the following way:
    \begin{enumerate}
        \item Set $k\coloneqq 0$, choose an initial approximation $f_0\in\Hilbert$.
        \item Find a dictionary element $d_{k+1}\in\Dictionary$ whose projection onto the subspace of $\Hilbert$, which is spanned by the current residual $r_k \coloneqq f-f_k$, is maximal, in other words \label{a:pga:2}
              \begin{align}
                  d_{k+1} & = \Argmax_{d\in\Dictionary} \Abs{\InnerProduct{r_k}{d}}.\label{e:pga:2}
              \end{align}
        \item Set
              \begin{align}
                  \alpha_{k+1} \coloneqq \InnerProduct{r_k}{d_{k+1}}.\label{e:pga:3}
              \end{align}
        \item Set $f_{k+1} \coloneqq f_k + \alpha_{k+1}d_{k+1}$.
        \item If a suitable stopping criterion is fulfilled: stop.\\
              Otherwise: increase $k$ by 1 and return to step~\ref{a:pga:2}.
    \end{enumerate}
    The approximation after $K$ steps is consequently given as $f_K = f_0 + \sum_{k=1}^K \alpha_k d_k$.
\end{algo}

Note that we use $\Hilbert=\Lebesgue{2}{\Sphere{d-1}}$ from now on.

The difficulty in the implementation of greedy algorithms in the application considered here lies in the computation of the inner products $\InnerProduct{r_k}{d} = \InnerProduct{f}{d}-\InnerProduct{f_k}{d}$ for $d\in\Dictionary$ in \eqref{e:pga:2} and \eqref{e:pga:3}, especially in the inner product $\InnerProduct{f}{d}$, whereas the inner product $\InnerProduct{f_k}{d}$ can be easily computed if the inner products $\InnerProduct{d}{d'}$ and $\InnerProduct{f_0}{d}$ are known for all $d,d'\in\Dictionary$.
The reason for this difficulty is that $f$ is unknown and that neither values of $f$ nor values of functionals evaluated at $f$ are given, which could be handled by a matching pursuit and a regularized functional matching pursuit, respectively.
Instead, realizations of the random variable $X$ are given, whose distribution possesses the density $f$.

To overcome this difficulty, we observe that by the definition of the inner product and the expectation value of a random variable, we have
\begin{align}
    \InnerProduct{f}{d} & = \int_{\Sphere{d-1}} d\B{\xi} f\B{\xi} \Differential{\SurfaceMeas}\B{\xi} = \Expectation{f}{d\B{X}},
\end{align}
since the distribution of $X$ possesses the density $f$.
Moreover, we use the strong law of large numbers (see, \eg, \cite[Chapter~5.3]{Klenke:2014}) to obtain the approximation
\begin{align}
    \InnerProduct{f}{d} & = \Expectation{f}{d\B{X}} \approx \frac{1}{N} \sum_{n=1}^N d\B{X_n},\label{e:slln}
\end{align}
which we insert into \eqref{e:pga:2} and \eqref{e:pga:3} to obtain the following algorithm.

\begin{algo}[Greedy algorithm for the estimation of PDFs]\label{a:pdfga}
    Let $f\in\Lebesgue{2}{\Sphere{d-1}}$ be an unknown PDF and $X_1,\ldots,X_N\in\Sphere{d-1}$ be a realization of a random variable $X$, whose distribution possesses the density $f$.
    Let a dictionary $\Dictionary$ be given as in Algorithm~\ref{a:pga}.

    Generate a sequence $\B{f_k}_{k\in\Natural_0}$ of approximations of $f$ iteratively according to the following scheme:
    \begin{enumerate}
        \item Set $k\coloneqq0$, choose an initial approximation $f_0\in\Lebesgue{2}{\Sphere{d-1}}$.
        \item Find a dictionary element $d_{k+1}\in\Dictionary$ fulfilling the maximization property\label{a:pdfga:2}
              \begin{align}
                  d_{k+1} & = \Argmax_{d\in\Dictionary} \Abs{\frac{1}{N}\sum_{n=1}^N d\B{X_n} - \InnerProduct{f_k}{d}}.\label{e:pdfga:2}
              \end{align}
        \item Set
              \begin{align}
                  \alpha_{k+1} & \coloneqq \frac{1}{N}\sum_{n=1}^N d_{k+1}\B{X_n} - \InnerProduct{f_k}{d_{k+1}}.\label{e:pdfga:3}
              \end{align}
        \item Set $f_{k+1} \coloneqq f_k + \alpha_{k+1}d_{k+1}$.
        \item If a suitable stopping criterion is fulfilled: stop.\\
              Otherwise: increase $k$ by 1 and return to step~\ref{a:pdfga:2}.
    \end{enumerate}
\end{algo}

Note that one would expect, from the theoretical point of view, a condition on the dictionary like $\Closure{\Span\Dictionary} = \Hilbert$.
However, from the computational point of view this condition can never be fulfilled, unless $\Hilbert$ is finite-dimensional, since $\Dictionary$ would have to be comprised of infinitely many elements, which cannot be realized on a computer.
The maximization in \eqref{e:pdfga:2} can consequently be accomplished by evaluating the term that has to be maximized for every dictionary element and choosing the dictionary element with the maximal value.
Note also that the maximizer in \eqref{e:pdfga:2} does not need to be unique.
In this case, we choose one arbitrary maximizer among all available maximizers.

By our motivation, the dictionary will consist of the kernels arising in the KDE and is consequently given as
\begin{align}
    \Dictionary & \coloneqq \Set{ \xi\mapsto K\B{X_n\cdot\xi} | n=1,\ldots,N }
\end{align}
for a given kernel $K$.

In this case, both in \eqref{e:pdfga:2} and \eqref{e:pdfga:3}, the kernel has to be evaluated at least $\B{N+1}N/2$ times, which---at the first sight---brings with it the same difficulties already discussed in Section~\ref{s:kde} if $N$ is large.
Two arguments can be given which reduce this difficulty in practice:

First, the greedy algorithm has to be applied only once per data set to obtain a sparse estimation of the PDF.
The coefficients $\alpha_k$ and the dictionary elements chosen by the greedy algorithm can be stored and reused, when an estimation of the PDF is needed in other algorithms.
On the contrary, when using a KDE itself, the kernel has to be evaluated very often inside these other algorithms.
This fact makes it possible to execute the greedy algorithm on a powerful and fast computer system, for example a large cluster, to obtain a list of chosen coefficients and kernels, and subsequently use a slower computer, for example an ordinary desktop computer, for the further applications.

Secondly, the term
\begin{align}
    \frac{1}{N}\sum_{n=1}^N d\B{X_n}\label{e:empexp}
\end{align}
arising in \eqref{e:pdfga:2} for every $d\in\Dictionary$,
which represents an empirical expectation value,
does not change in the course of the iteration.
That is to say, \eqref{e:empexp} has to be computed and stored only once in the beginning of the algorithm and does not need to be recomputed in every iteration.
The value computed prior to the iteration can then also be used in \eqref{e:pdfga:3}.
Additionally, the value of \eqref{e:empexp} for some $d$ can be computed independently from the other dictionary elements.
Thus, this process is ideal for parallelization, which improves the performance on those fast computer systems mentioned before.

Note that, in principle, this method can be extended to other types of dictionary elements apart from kernels, for example, spherical harmonics or wavelets, which have not been implemented for the numerical experiments in this paper.
It would also be possible to use kernels with different concentration parameters for the dictionary.
The use of such overcomplete dictionaries turned out to be useful for spherical inverse problems in the geosciences (see \cite{Berkel:2011,Fischer:2011,Fischer:2012,Fischer:2013a,Fischer:2013,Michel:2015,Michel:2014,Michel:2016,Telschow:2014}).
Further note that greedy algorithms with dictionaries consisting of kernels also appear in a method called \emph{kernel matching pursuit} in the field of machine learning (\cf{} \cite{Vincent:2002}).

\section{Numerical experiments}\label{s:num}
The Abel-Poisson kernel as introduced in \eqref{e:ap:def} was chosen with a parameter $h=0.9$ for all of the following numerical experiments.
Since Algorithm~\ref{a:pdfga} uses the approximation in \eqref{e:slln}, and convergence for the PGA is only known if the inner product can be evaluated exactly, a synthetic example is presented first in this section to empirically study the convergence of the algorithm.
We proceed by a more detailed description of the CT data set used for the following computations.
We continue by presenting simulation results using the greedy approximation of the PDF and comparing computation times between the use of KDEs and the greedy approximation.
Finally, we conclude this section by presenting some first results about the validation of fiber laydown models.
For this purpose, we also use the greedy algorithm.

\subsection{A synthetic example}
To create a synthetic example to test the convergence of the greedy algorithm, we choose the PDF $f\B{\xi} \coloneqq \AbelPoissonK{3}{0.6}\B{\xi\cdot\eps^3}$, where the Abel-Poisson kernel is concentrated around $\eps^3 = \B{0,0,1}^\Transpose$ and apply the acceptance-rejection method (\cf{} section~\ref{s:kde} of this paper) to sample a data set $X_1,\ldots,X_N$ from this PDF, where $N=10^6$.
The PDF and \num{2000} points from the data set are presented in Figure~\ref{f:syn:pdf,pts}.

\begin{figure}
    \centering
    \input{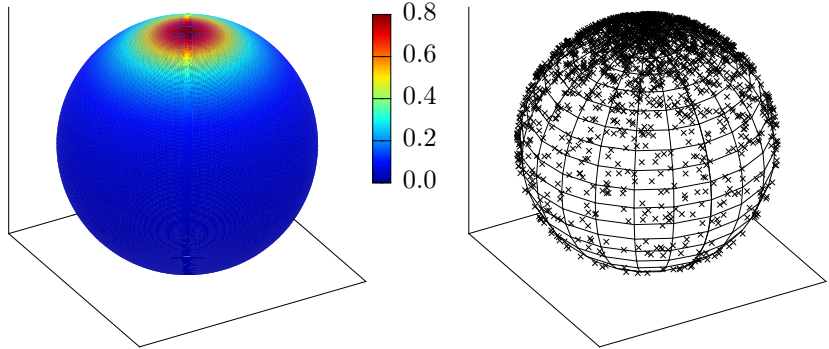}
    \caption{Left: PDF from the synthetic example, an Abel-Poisson kernel with parameter $h=0.6$, centered at $\eps^3=\B{0,0,1}^\Transpose$;
        Right: 2000 of in total $10^6$ data points sampled from that PDF with the acceptance-rejection method.}\label{f:syn:pdf,pts}
\end{figure}

\begin{figure}
    \centering
    \input{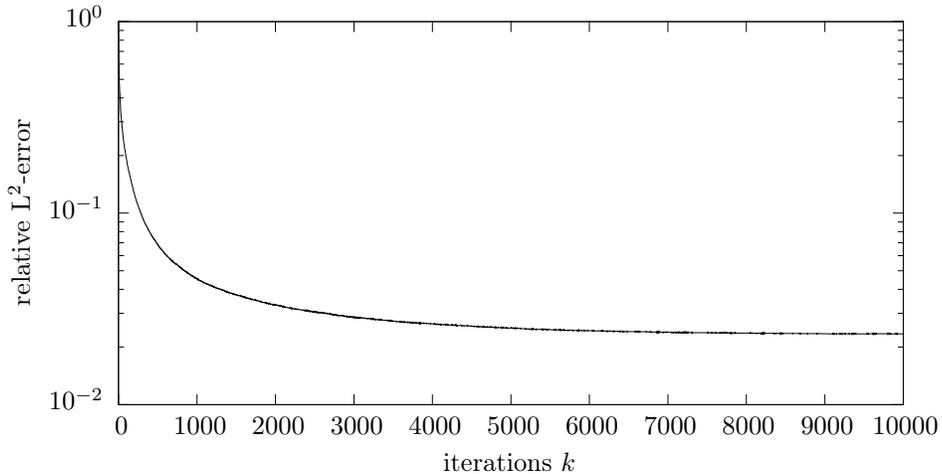}
    \caption{Semi-logarithmic plot of the relative $\Lebesgue{2}{\Sphere{d-1}}$-error $\Norm{f-f_k}/\Norm{f}$ of the first \num{10000} iterations of the greedy algorithm, when applied to the synthetic example.}\label{f:syn:err}
\end{figure}

In this synthetic case, where the PDF is already known, it is possible to compute the relative $\Lebesgue{2}{\Sphere{d-1}}$-error $\Norm{f-f_k}/\Norm{f}$ explicitly due to properties of the Abel-Poisson kernel.
A semi-logarithmic plot of this relative error for \num{10000} iterations of the greedy algorithm can be found in Figure~\ref{f:syn:err}.

The algorithm shows a convergent behavior, at least in these first \num{10000} iterations.
Multiple performed experiments show that the number of data points influences the convergence of the algorithm in the following way:
the more iterations are performed, the more data points are needed such that the algorithm shows a convergent behavior.
This is plausible, since, having reached a certain level of accuracy, the approximation can only be improved with additional information on the function to be approximated.
The dependence of the convergence on the number of data points is a subject of future research.

\subsection{Description of the CT data set}\label{s:data}
We deal with a spunbond process of the industrial company Oerlikon Neumag.
After specifying several production parameters, the company produced samples of nonwovens which will be used for the presented numerical experiments.
The samples were analyzed by the department \emph{Image Processing} of the Fraunhofer ITWM with its 3D-microtomography scanner.
The CT scan of the real samples delivers a real-valued third-order tensor with so-called gray values. 
The gray values are mapped to local fiber orientations at each voxel with the help of an eigenvalue analysis of the Hessian matrix of the second partial derivatives of the gray values.
For more details, the reader is referred to \cite{Redenbach:2012}.
As already mentioned in the introduction, the microtomography scanner possesses a resolution of \SIrange{1}{70}{\micro\metre} while the sample size is \SIrange{1}{10}{\cubic\milli\metre}.
For the data set that we use, this high resolution leads to $N=\num{9600558}$ points on $\Sphere{2}$, which will be the basis for our considerations.
\subsection{Application of the greedy and simulation algorithms}
\begin{figure}
        \centering
        \input{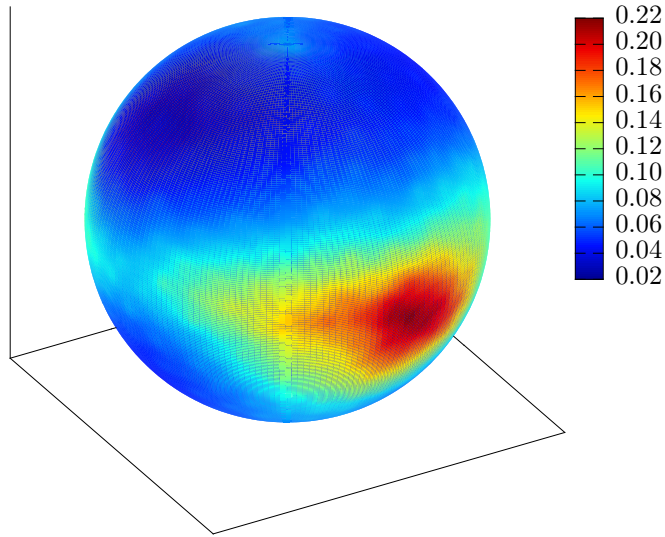}
        \caption{Sparse estimation of the PDF of a real data set after \num{10000} iterations of the greedy algorithm.}\label{f:ct:pdfga}
\end{figure}

We applied Algorithm~\ref{a:pdfga} on the CT data set which was described in the previous section.
The result is a sparse kernel density estimator, which is depicted in Figure~\ref{f:ct:pdfga}.

\begin{figure}
    \centering
    \input{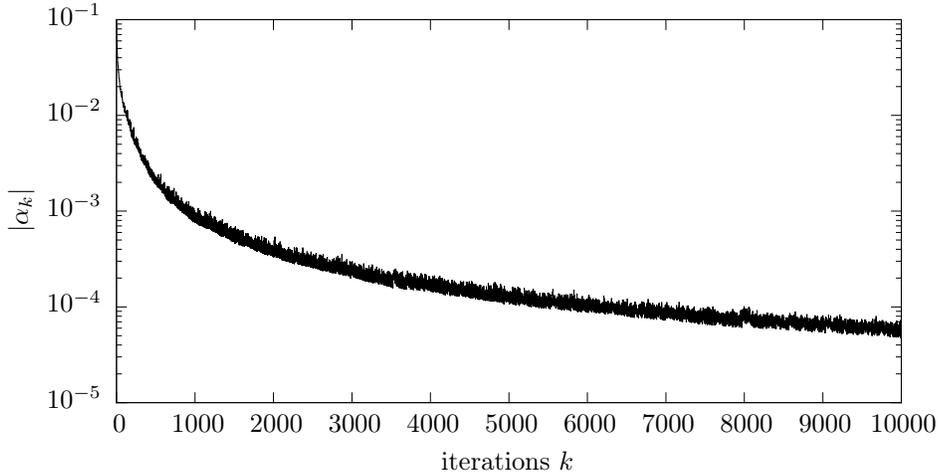}
    \caption{Absolute values of the chosen coefficients in the first \num{10000} iterations of the greedy algorithm, applied to the real data set.}\label{f:ct:err}
\end{figure}

Since the error of the estimation can no longer be explicitly computed, in Figure~\ref{f:ct:err} the absolute values of the chosen coefficients $\Abs{\alpha_k}$ are plotted, to get an impression of the size of the corrections which the greedy algorithm performs in each of the first \num{10000} iterations.
This indicator for the size of corrections also shows a convergent behavior, similarly to the synthetic example above.

Note that the KDE in Figure~\ref{f:ct:kde} and the sparse estimator in Figure~\ref{f:ct:pdfga} show a clear qualitative similarity.
Quantitatively, one recognizes certain differences both in the values of the estimator (the maximum is \num{0.18} for the KDE and \num{0.22} for the greedy approximation) and in the structure of the estimator (the structures are rougher in the result of the greedy approximation).
This is no contradiction, since the greedy approximation and the KDE are both estimations of the PDF, but the greedy approximation is no approximation of the KDE.
Thus, depending on the choice of parameters $h$ both for the KDE and the dictionary of the greedy algorithm, both estimations can be different while still estimating the same PDF.

\begin{figure}
    \centering
    \input{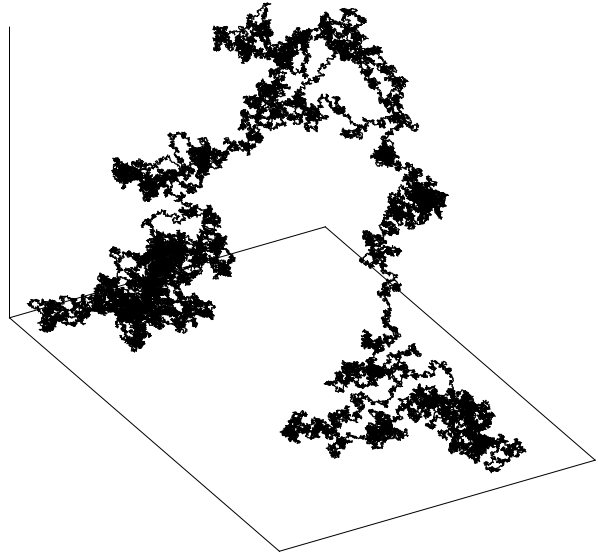}
    \caption{Fiber, consisting of \num{100000} points, simulated by Algorithm~\ref{a:sim} in conjunction with a greedy approximation of the PDF.}\label{f:sim}
\end{figure}

Remember that our initial motivation for the development of a greedy algorithm was the inefficiency of Algorithm~\ref{a:sim} if the PDF is estimated by a kernel density estimator.
For this reason, we use the approximation $f_{\num{10000}}$ (generated by the greedy algorithm above) for a 3D simulation of nonwoven fabrics.
A discretized fiber that has been generated this way is shown in Figure~\ref{f:sim}.

A comparison of computation times of the simulation algorithm is shown in Table~\ref{t:times}.

In this table, we compare the uses of the KDE and the sparse KDE as well as the application of two different upper bounds for the PDF in the acceptance-rejection method (see the considerations in section~\ref{s:kde}).
The following raw data are listed in the table: the number of samples $X_j$ that have been generated, the number of evaluations of the estimated PDF that were needed to generate these samples, and the CPU time that was consumed.
Note that there is a difference of the number of samples between the four presented scenarios due to the high computational effort that would be needed for the more inefficient methods.
For that reason, to achieve a better comparison, we add the average number of evaluations of the estimated PDF per sample, the average CPU time per sample, and the time that would be needed to simulate a nonwoven fabric with \num{100} fibers with \num{100000} segments each, to the table.
\begin{table}\centering
    \caption{Computation times for different estimates of the PDF and different upper bounds for the estimated PDF.
        For the comparison \enquote{CPU time/nonwoven}, it has been assumed that a simulated nonwoven consists of \num{100} fibers with \num{100000} segments each. The presented data show an enormous saving of computation time when using the newly developed greedy algorithm.}
    \label{t:times}
    \begin{tabular}{lrrrr}\toprule
        & \multicolumn{2}{c}{KDE} & \multicolumn{2}{c}{Greedy approximation}\\
        upper bound & triangle ineq. & evaluation & triangle ineq. & evaluation\\\midrule
        samples & \num{420} & \num{25000} & \num{100000} & \num{100000}\\
        evaluations & \num{85607} & \num{59894} & \num{258022} & \num{166128} \\
        CPU time    & \SI{16022}{\second} & \SI{8637.3}{\second} & \SI{37.543}{\second} & \SI{24.903}{\second}\\\midrule
        eval./sample & \num{204} & \num{2.40} & \num{2.58} & \num{1.66}\\
        CPU time/sample & \SI{3.81e1}{\second} & \SI{3.45e-1}{\second} & \SI{3.75e-4}{\second} & \SI{2.49e-4}{\second}\\
        CPU time/nonwoven & \SI{105833}{\hour} & \SI{958}{\hour} & \SI{1.04}{\hour} & \SI{0.69}{\hour}\\\bottomrule
    \end{tabular}
\end{table}
One can see that the use of the greedy approximation makes Algorithm~\ref{a:sim} more efficient by multiple orders of magnitude.
Moreover, it is clear that the upper bound that is determined by applying the triangle inequality is much worse for the KDE, since there are much more summands in the estimator, such that a comparison of these numbers is not fair.
Nevertheless, if the more efficient variant by evaluation on a fine grid is used for both estimators, the results are comparable.
In that case, the simulation time for a nonwoven with \num{100} fibers with \num{100000} line segments each, is reduced from \SI{958}{\hour} to \SI{0.69}{\hour}, thus from nearly \num{40} days to \num{41} minutes, a factor of nearly \num{1400}.

Unfortunately, as in the case of the other known simulation methods for nonwovens (\cf{} the introduction), it is not trivial to validate the output of the simulation algorithm.
By construction of the method, the output will fulfill the validation criterion introduced in the following section.
Nevertheless, further research should be conducted on the topic of evaluating the quality of the simulations.

\subsection{Validation of fiber laydown models}\label{s:val}
Another problem that can be attacked by using the developed greedy algorithm is the validation of a fiber laydown model, which has been developed at Fraunhofer ITWM and which is implemented in their simulation tool FIDYST (see \cite{Gramsch:2015}).
For this purpose, we analyze the same data set generated by CT scans of real nonwoven fabrics as above and the result of a FIDYST simulation using the greedy algorithm introduced before.

As already mentioned before, the company Oerlikon Neumag has specified certain production parameters, which were used for the fabrication of the samples.
For example, the belt speed was set to \SI{38}{\meter\per\minute}, the mean spinning speed is \SI{79.4}{\meter\per\second}, and the fibers are made of polypropylene with a fiber diameter of \SI{1.2e-5}{\meter}.
All specified parameters are used as input to FIDYST, which simulates the air flow and the fiber dynamics of this process---until the filaments are laid down on the conveyor belt.
We prescribed a total simulation time of \SI{0.5}{\second}.
Due to the mean spinning speed, we end up with a final filament length of \SI{39.7}{\meter}.
For the upcoming analysis, we neglect the fiber parts in the air and restrict ourselves to the laid-down parts on the conveyor belt. 
With a spatial discretization of the fibers of \SI{2.5e-4}{\meter} this leads to approximately \num{150000} points per simulated fiber.

Since FIDYST can only produce a two-dimensional output, we face the problem to compare a 2D data set from the simulation with a real 3D data generated by CT scans.
Let $\mathbb{X}=\set{X_1,\ldots,X_N}\subset\Sphere{2}$ denote the CT data set and $\mathbb{Y}=\set{Y_1,\ldots,Y_M}\subset\Sphere{1}$ the FIDYST data set, where $N,M\in\Natural$ are the numbers of data points for both types of data.
It is known from the measurement process, that the third component of the CT data set corresponds to the direction perpendicular to the moving belt in the production process.
Consequently, we overcome the challenge of comparing data sets in different dimensions by generating a 2D data set $\tilde{\mathbb{X}}= \set{\tilde{X}_1,\ldots,\tilde{X}_N}\subset\Sphere{1}$ by projection onto the first two components.

\begin{figure}
    \input{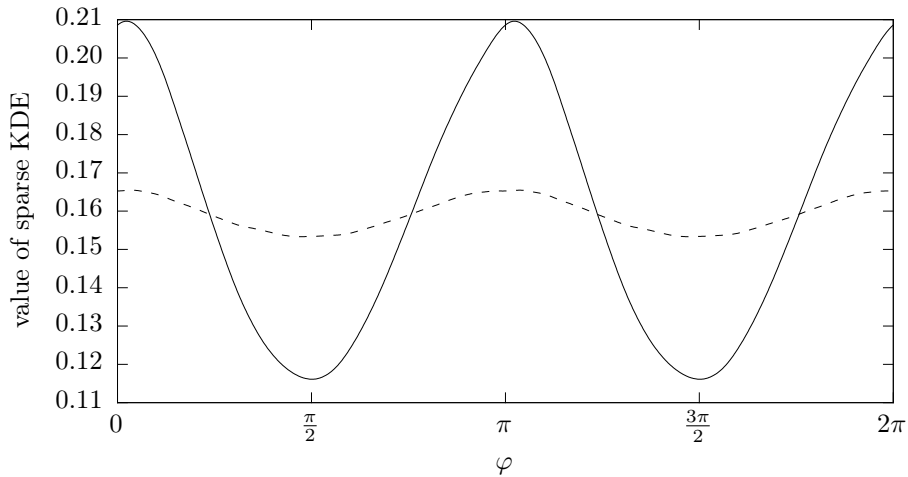}
    \caption{Sparse KDEs of the CT data set $\tilde{\mathbb{X}}$ (solid line) and the FIDYST data set $\mathbb{Y}$ (dashed line).
        Two-dimensional Abel-Poisson kernels with $h=0.9$ were used in both cases.
        The sphere $\Sphere{1}$ is identified with the interval $[0,2\pi)$ on the horizontal axis.}
        \label{f:val1}
\end{figure}

The data sets used in the following consist of $N=\num{9600558}$ (CT data) and $M=\num{3137432}$ (FIDYST data) points, where we concatenated data from \num{10} distinct simulations (all initialized with the same parameters), and symmetrized the data.
The application of the algorithm on $\tilde{\mathbb{X}}$ and $\mathbb{Y}$ yields the results depicted in Figure~\ref{f:val1}.
A possible interpretation of Figure~\ref{f:val1} is as follows:
it can be seen that the distribution of filament directions of the simulated fiber possesses maxima and minima at the same angles as the real fibers.
However, the variance of the directions is much larger in the CT data set than in the FIDYST data set.
Multiple calculations show that this effect also arises with other process parameters.
It seems to be the case that the model used by FIDYST could be improved to better match the real world data.
Necessary modifications of the model are the subject of further research.

Note that, so far, there was no way to measure the consistency of the model with reality.
For the first time, the introduced method offers a possibility to do so.
\section{Conclusions} 
In this paper, we dealt with the three-dimensional simulation of nonwovens.
Using data from computer tomography scans of real nonwoven fabrics, we proposed an algorithm for the simulation of nonwovens.
The basis of this algorithm is an estimate of the distribution of fiber directions inside the nonwoven, that is a probability density function.
It turned out that using kernel density estimators, which are a standard technique in nonparametric statistics, leads to a very high computational effort when applying the simulation algorithm due to the high amount of available data.
By introducing a greedy algorithm, which yields a sparse estimator of the unknown probability density, the computation time of the fiber simulation algorithm could be reduced by a factor of \num{1400}, from more than a month to less than an hour.

The results presented in the work at hand give rise to several other topics that should be addressed in future research.
First, the convergence of the greedy algorithm should be analyzed theoretically, in particular the dependence of the convergence on the number of data points.
Secondly, in this paper, we restricted ourselves to a dictionary consisting of Abel-Poisson kernels centered at the data points in analogy to kernel density estimators.
Of course, other functions could be used, for example Abel-Poisson kernels with different parameters $h$, that may also be centered at a regular grid instead of the data points.
We could also add global functions like spherical harmonics to the dictionary.
Thirdly, the presented simulation algorithm is not limited to the spunbond process.
Several other production processes for nonwovens or also other technical textiles could be simulated with the given method, as long as directional CT data are available.
Additionally, the simulation algorithm could be extended to include different belt movements than the one, which is already incorporated in the data set.
Finally, the approach for the validation of fiber laydown models should be pursued further.
The reason for the inconsistency of the distribution of fiber directions between the FIDYST model and the CT data should be analyzed.
Also other fiber laydown models, like the stochastic surrogate model, could be analyzed with the given validation approach.
Furthermore, different models, like the FIDYST and the surrogate model, could be compared to each other, instead to the CT data, by our approach.

\bibliographystyle{amsplain}
\bibliography{nonwoven.bib}

\providecommand{\bysame}{\leavevmode\hbox to3em{\hrulefill}\thinspace}
\providecommand{\MR}{\relax\ifhmode\unskip\space\fi MR }
\providecommand{\MRhref}[2]{%
  \href{http://www.ams.org/mathscinet-getitem?mr=#1}{#2}
}
\providecommand{\href}[2]{#2}
\begin{thebibliography}{10}

\bibitem{Albrecht:2003}
W~Albrecht, H~Fuchs, and W~Kittelmann (eds.), \emph{{Nonwoven Fabrics: Raw
  Materials, Manufacture, Applications, Characteristics, Testing Processes}},
  Wiley-VCH, Weinheim, 2003.

\bibitem{Banerjee:2010}
A~Banerjee and P~Burlina, \emph{Efficient particle filtering via sparse kernel
  density estimation}, IEEE Trans. Image Process. \textbf{19} (2010),
  2480--2490.

\bibitem{Barron:2008}
AR~Barron, A~Cohen, W~Dahmen, and RA~DeVore, \emph{Approximation and learning
  by greedy algorithms}, Ann.\ Statist. \textbf{36} (2008), 64--94.

\bibitem{Berkel:2011}
P~Berkel, D~Fischer, and V~Michel, \emph{Spline multiresolution and numerical
  results for joint gravitation and normal-mode inversion with an outlook on
  sparse regularisation}, Int. J. Geomath. \textbf{1} (2011), 167--204.

\bibitem{Cacoullos:1966}
T~Cacoullos, \emph{Estimation of a multivariate density}, Ann. Inst. Statist.
  Math. \textbf{18} (1966), 179--189.

\bibitem{Chen:2013}
F~Chen, H~Yu, J~Yao, and R~Hu, \emph{Robust sparse kernel density estimation by
  inducing randomness}, Pattern Anal. Appl. \textbf{18} (2013), 367--375.

\bibitem{Devroye:1986}
L~Devroye, \emph{{Non-Uniform Random Variate Generation}}, Springer, New York,
  1986.

\bibitem{Fischer:2011}
D~Fischer, \emph{{Sparse Regularization of a Joint Inversion of Gravitational
  Data and Normal Mode Anomalies}}, Ph.D. thesis, Geomathematics Group,
  University of Siegen, 2011, published by Dr.\ Hut, München.

\bibitem{Fischer:2012}
D~Fischer and V~Michel, \emph{Sparse regularization of inverse
  gravimetry---case study: spatial and temporal mass variations in {South
  America}}, Inverse Problems \textbf{28} (2012), 065012.

\bibitem{Fischer:2013a}
\bysame, \emph{Automatic best-basis selection for geophysical tomographic
  inverse problems}, Geophys. J. Int. \textbf{193} (2013), 1291--1299.

\bibitem{Fischer:2013}
\bysame, \emph{Inverting {GRACE} gravity data for local climate effects}, J.
  Geod. Sci. \textbf{3} (2013), 151--162.

\bibitem{Foucart:2013}
S~Foucart and H~Rauhut, \emph{{A Mathematical Introduction to Compressive
  Sensing}}, Birkhäuser, New York, 2013.

\bibitem{Freeden:1981}
W~Freeden, \emph{On spherical spline interpolation and approximation}, Math.
  Method. Appl. Sci. \textbf{3} (1981), 551--575.

\bibitem{Freeden:2002}
W~Freeden and K~Hesse, \emph{On the multiscale solution of satellite problems
  by use of locally supported kernel functions corresponding to equidistributed
  data on spherical orbits}, Studia Sci. Math. Hungar. \textbf{39} (2002),
  37--74.

\bibitem{Freeden:2009}
W~Freeden and M~Schreiner, \emph{{Spherical Functions of Mathematical
  Geosciences: A Scalar, Vectorial, and Tensorial Setup}}, Springer, Berlin,
  2009.

\bibitem{Friedman:1981}
JH~Friedman and W~Stuetzle, \emph{Projection pursuit regression}, J. Amer.
  Statist. Assoc. \textbf{76} (1981), 817--823.

\bibitem{Friedman:1984}
JH~Friedman, W~Stuetzle, and A~Schroeder, \emph{Projection pursuit density
  estimation}, J. Amer. Statist. Assoc. \textbf{79} (1984), 599--608.

\bibitem{Galassi:2009}
M~Galassi, J~Davies, J~Theiler, B~Gough, G~Jungman, P~Alken, M~Booth, and
  F~Rossi, \emph{{GNU Scientific Library Reference Manual}}, 3 ed., Network
  Theory, Bristol, 2009.

\bibitem{Gramsch:2015}
S~Gramsch, D~Hietel, and R~Wegener, \emph{Optimizing spunbond, meltblown, and
  airlay processes with {FIDYST}}, Melliand Int. \textbf{21} (2015), 115--117.

\bibitem{Grothaus:2014}
M~Grothaus, A~Klar, J~Maringer, P~Stilgenbauer, and R~Wegener,
  \emph{Application of a three-dimensional fiber lay-down model to non-woven
  production processes}, J. Math. Ind. \textbf{4{\normalfont{}:4}} (2014).

\bibitem{Goetz:2007}
T~Götz, A~Klar, N~Marheineke, and R~Wegener, \emph{A stochastic model and
  associated {Fokker-Planck} equation for the fiber lay-down process in
  nonwoven production processes}, SIAM J. Appl. Math. \textbf{67} (2007),
  1704--1717.

\bibitem{Hall:1987}
P~Hall, GS~Watson, and J~Cabrera, \emph{Kernel density estimation with
  spherical data}, Biometrika \textbf{74} (1987), 751--762.

\bibitem{Hong:2010}
X~Hong, S~Chen, and CJ~Harris, \emph{Sparse kernel density estimation technique
  based on zero-norm constraint}, {The 2010 International Joint Conference on
  Neural Networks}, IEEE, 2010, IJCNN, Barcelona, 18--23 July 2010.

\bibitem{Jones:1996}
MC~Jones, JS~Marron, and SJ~Sheather, \emph{A brief survey of bandwidth
  selection for density estimation}, J. Amer. Statist. Assoc. \textbf{91}
  (1996), 401--407.

\bibitem{Klar:2009}
A~Klar, N~Marheineke, and R~Wegener, \emph{Hierarchy of mathematical models for
  production processes of technical textiles}, ZAMM Z. Angew. Math. Mech.
  \textbf{89} (2009), 941--961.

\bibitem{Klar:2012}
A~Klar, J~Maringer, and R~Wegener, \emph{A 3d model for fiber lay-down in
  nonwoven production processes}, Math. Models Methods Appl. Sci. \textbf{22}
  (2012), 1250020.

\bibitem{Klenke:2014}
A~Klenke, \emph{{Probability Theory}}, Springer, London, 2014.

\bibitem{Kloeden:1992}
P~Kloeden and E~Platen, \emph{{Numerical Solution of Stochastic Differential
  Equations}}, Springer, Berlin, 1992.

\bibitem{Mallat:1993}
SG~Mallat and Z~Zhang, \emph{Matching pursuits with time-frequency
  dictionaries}, IEEE Trans. Signal Process. \textbf{41} (1993), 3397--3415.

\bibitem{Mao:2010}
N~Mao and SJ~Russell, \emph{Modelling of nonwoven materials}, Modelling and
  predicting textile behaviour (X~Chen, ed.), Woodhead, Cambridge, 2010,
  pp.~180--224.

\bibitem{Marheineke:2011}
N~Marheineke and R~Wegener, \emph{Modeling and application of a stochastic drag
  for fiber dynamics in turbulent flows}, Int. J. Multiphase Flow \textbf{37}
  (2011), 136--148.

\bibitem{Michel:2013}
V~Michel, \emph{{Lectures on Constructive Approximation}}, Birkhäuser, Basel,
  2013.

\bibitem{Michel:2015}
\bysame, \emph{{RFMP}: An iterative best basis algorithm for inverse problems
  in the geosciences}, {Handbook of Geomathematics} (W~Freeden, MZ~Nashed, and
  T~Sonar, eds.), Springer, Berlin, 2 ed., 2015, pp.~2121--2147.

\bibitem{Michel:2014}
V~Michel and R~Telschow, \emph{A non-linear approximation method on the
  sphere}, Int. J. Geomath. \textbf{5} (2014), 195--224.

\bibitem{Michel:2016}
\bysame, \emph{The regularized orthogonal functional matching pursuit for
  ill-posed inverse problems}, SIAM J. Numer. Anal. \textbf{54} (2016),
  262--287.

\bibitem{Parzen:1962}
E~Parzen, \emph{On estimation of a probability density function and mode}, Ann.
  Math. Stat. \textbf{33} (1962), 1065--1076.

\bibitem{Redenbach:2012}
C~Redenbach, A~Rack, K~Schladitz, O~Wirjadi, and M~Godehardt, \emph{Beyond
  imaging: on the quantitative analysis of tomographic volume data}, Int. J.
  Mater. Res. \textbf{103} (2012), 217--227.

\bibitem{Rosenblatt:1956}
M~Rosenblatt, \emph{Remarks on some nonparametric estimates of a density
  function}, Ann. Math. Stat. \textbf{27} (1956), 832--837.

\bibitem{Rosset:2002}
S~Rosset and E~Segal, \emph{Boosting density estimation}, Advances in Neural
  Information Processing Systems 15 (S~Becker, S~Thrun, and K~Obermayer, eds.),
  MIT Press, 2003, NIPS, Vancouver, 9--14 December 2002, pp.~657--664.

\bibitem{Telschow:2014}
R~Telschow, \emph{{An Orthogonal Matching Pursuit for the Regularization of
  Spherical Inverse Problems}}, Ph.D. thesis, Geomathematics Group, University
  of Siegen, 2014, published by Dr.\ Hut, München.

\bibitem{Temlyakov:2008}
VN~Temlyakov, \emph{Greedy approximation}, Acta Numer. \textbf{17} (2008),
  235--409.

\bibitem{Verbeek:2003}
JJ~Verbeek, N~Vlassis, and B~Kröse, \emph{Efficient greedy learning of
  gaussian mixture models}, Neural Comput. \textbf{15} (2003), 469--485.

\bibitem{Vincent:2002}
P~Vincent and Y~Bengio, \emph{Kernel matching pursuit}, Mach. Learn.
  \textbf{48} (2002), 165--187.

\bibitem{Wegener:2015}
R~Wegener, N~Marheineke, and D~Hietel, \emph{Virtual production of filaments
  and fleeces}, {Currents in Industrial Mathematics} (H~Neunzert and
  D~Prätzel-Wolters, eds.), Springer, Berlin, 2015, pp.~103--162.

\end{thebibliography}
\end{document}